\documentclass{amsart}

\usepackage{amssymb}
 \usepackage{amsmath}
\usepackage{amsfonts}
\usepackage{graphicx} 

\theoremstyle{plain}
 \newtheorem{thm}{Theorem}
 \newtheorem{lemma}[thm]{Lemma}
 \newtheorem{prop}[thm]{Proposition}

\newtheorem{attempt}[thm]{Attempt}
\theoremstyle{remark}
 \newtheorem*{remark}{Remark}

\newcommand{\cA}{\mathcal{A}}

 \newcommand{\cH}{\mathcal{H}}
 
  \newcommand{\cM}{\mathcal{M}}
  \newcommand{\cP}{\mathcal{P}}
    \newcommand{\cS}{\mathcal{S}}

 \newcommand{\bN}{\mathbb{N}}
  \newcommand{\bZ}{\mathbb{Z}}
    \newcommand{\bM}{\mathbb{M}}

 \newcommand{\al}{\alpha}

 \newcommand{\ep}{\varepsilon}

 \newcommand{\hot}{\hat{\otimes}}
 \newcommand{\vot}{\Bar{\otimes}}

\def\wk*{\hbox{\rm wk*}}

\title[Type I von Neumann]{A note on approximate amenability of type I von Neumann algebras}

\date{March 15, 2025}

\author{Yong Zhang}
\address{ Department of Mathematics,
University of Manitoba, Winnipeg R3T 2N2, Canada}
\email{yong.zhang@umanitoba.ca}

\begin{document}

\begin{abstract}
Using the methods of Ozawa \cite{Oz} and Runde \cite{Runde}, we show that a type I von Neumann algebra is approximately amenable if and only if it is amenable.
\end{abstract}

\maketitle

A Banach algebra $\cA$ is called \emph{approximately amenable} if for each Banach $\cA$-bimodule $X$, every continuous derivation $D$ from $\cA$ into $X$ is approximately inner \cite{G-L04, GLZ08}. That is, there is a net $(x_i)\subset X$ such that
\[
D(a) = \lim_i a\cdot x_i -x_i\cdot a \quad (a\in \cA)
\]
in the norm topology of $X$.

If $\cA$ has a central bounded approximate identity, then it is approximately amenable if and only if it is \emph{pseudo-amenable} \cite{Z22}. The latter means that there exists a net $(u_\al)\subset \cA\hot^\gamma \cA$ such that
\[
au_\al - u_i a \to 0, \quad \varDelta_\cA(u_\al)a\to a
\]
for all $a\in \cA$. Here $\hot^\gamma$ represents the Banach space projective tensor product, $\varDelta_\cA$: $\cA\hot^\gamma \cA\to \cA$ is the product mapping defined by $\varDelta_\cA (a\otimes b)= ab$, and the convergences are all in the norm topology \cite{GZ07}.
In particular, if $\cA$ has an identity, then the approximate amenability is the same as the pseudo-amenability for $\cA$ (see also \cite{Runde}).

It is an open question whether an approximately amenable C*-algebra must be amenable. In particular, we don't know whether a von Neumann algebra must be amenable if it is approximately amenable. Ozawa \cite{Oz} showed that $B(\cH)$ is not approximately amenable for any Hilbert space $\cH$ of infinite dimension. Runde gives a more friendly proof to this result in \cite{Runde}. In this note, we employ the method developed in \cite{Oz,Runde} to examine type I von Neumann algebras. We show that a type I von Neumann algebra is approximately amenable if and only if it is amenable.

Recall that a von Neumann algebra $\cM$ is of \emph{type I} if every non-trivial central projection $P\in \cM$ majorizes a non-trivial abelian projection. It is well-known that a von Neumann algebra $\cM$ is of type I if and only if it is in the form
\begin{equation}\label{E: type I}
    \cM = \oplus_{\al\in \Gamma} A_\al\vot B(\cH_\al),
\end{equation}

where $\oplus$ stands for an $\ell^\infty$ sum and, for each $\al$, $A_\al$ is a commutative von Neumann algebra and $\cH_\al$ is a Hilbert space. (See \cite[Theorem~V.1.27]{Tak I}.)

\begin{prop}\label{P: dim bdd}
Let $\cM$ be a type I von Neumann algebra of the form \eqref{E: type I}.
If there is $N\in \bN$ such that dim$(\cH_\al) \leq N$ for all $\al$, then $\cM$ is amenable.
\end{prop}

\proof

Since dim$\cH_\al\leq N$ ($\al\in \Gamma$), each $\cH_\al$ is equal to some $\cH_n$ with $n\leq N$. This implies
\[
\cM = \oplus_{k=1}^N A_k\otimes B(\cH_k)
\]
with $A_k$ ($1\leq k\leq N$) being a commutative von Neumann algebra.
Therefore, $\cM$ is amenable as a finite sum of amenable algebras.
\qed

Let $\cH$ be a Hilbert space and $(\cH_\al)_{\al\in \Gamma}$ be a collection of orthogonal subspaces of $\cH$ such that 
\[
\cH =\ell^2\text{-}\oplus_{\al\in \Gamma}\cH_\al .
\]
Then, with a natural embedding, we may regard $\ell^\infty$- $\oplus_{\al\in \Gamma}B(\cH_\al)$ as a subspace of $B(\cH)$. For an element $a=\oplus_{\al\in \Gamma}a_\al \in \oplus_{\al\in \Gamma} B(\cH_\al)$, we have $|a| = \oplus_{\al\in \Gamma}|a_\al|$.

Denote by $L_2(\cH)$ and $L_1(\cH)$ the space of Hilbert-Schmidt operators and the space of trace class operators on $\cH$, respectively. The Hilbert-Schmidt norm and the trace class norm of $a\in B(\cH)$ are denoted by $\|a\|_2$ and $\|a\|_1$, respectively. Let $a=\oplus_{\al\in \Gamma}a_\al \in \oplus_{\al\in \Gamma} B(\cH_\al)$.
It is readily checked that
\[
\|a\|_2 = (\sum_{\al\in \Gamma}\|a_\al\|_2^2)^\frac{1}{2}
\]
and 
\[
\|a\|_1 = \sum_{\al\in \Gamma}\|a_\al\|_1.
\]
We then have
\[
\oplus_{\al\in \Gamma} B(\cH_\al) \cap L_2(\cH) = \ell^2 \text{-} \oplus_{\al\in \Gamma} L_2(\cH_\al)
\]
and
\[
\oplus_{\al\in \Gamma} B(\cH_\al) \cap L_1(\cH) = \ell^1 \text{-} \oplus_{\al\in \Gamma} L_1(\cH_\al).
\]

Let $E_\al$ be an orthonormal basis of $\cH_\al$. Then $E = \cup_{\al\in \Gamma} E_\al$ is an orthonormal basis of $\cH$. As is well known $L_2(\cH) = \ell^2(E)\otimes_2 \ell^2(E)$ and $L_1(\cH) = \ell^2(E)\hot^\gamma \ell^2(E)$, where the former is the Hilbert space tensor product, the latter is the Banach space projective tensor product of $\ell^2(E)$ with itself. Let $F$ be a finite subset of $E$ and $P_F\in B(\cH)$ be the orthogonal projection from $\cH$ onto the subspace $\cH_F = \ell^2(F)$. With this setting, \cite[Lemma~2.1]{Oz} for the case $p=2$ leads to the following lemma.
\begin{lemma}\label{L: Oz}
    Let $T=\sum_{i=1}^ra_i\otimes b_i\in B(\cH)\otimes B(\cH)$.  Then for each $e\in E$, 
    \[
   T_F(e) = \sum_{i=1}^rP_F \cdot a_i(e)\otimes P_F\cdot b_i^*(e)\in \ell^2(F)\hot^\gamma \ell^2(F)
    \]
    and
    \[
    \sum_{e\in E} \|T_F(e)\|_{\ell^2(F)\hot^\gamma \ell^2(F)} \leq |F| \|T\|_\gamma,
    \]
    where $\|T\|_\gamma$ is the projective tensor norm of $T$ as an element of $B(\cH)\hot^\gamma B(\cH)$ .
    \end{lemma}

Let $G= SL_3(\bZ)$. As is well known, $G$ is a finitely generated group with Kazhdan's property (T). 
Let $(p_k)_{k=1}^\infty$ be the increasing sequence of all prime numbers, and let $\cP_k$ be the projective plane over the field $\bZ/p_k\bZ$ obtained from $(\bZ/p_k\bZ)^3$.
We have $|\cP_k| = p_k(p_k + 1) +1$. Following the construction of Ozawa, we see that
the natural unitary representation of $G$ on $\bZ^3$ induces an invertible action $\pi$ of $G$ on $\cP_k$. (For each $x\in G$, $\pi(x)$ is indeed a permutation of $\cP_k$.) This, in turn, induces a unitary representation of $G$ on $\ell^2(\cP_k)$. We denote it by $\pi_k$. Let $x_1,x_2,\ldots, x_m$ be fixed generators of $G$. 
We have $\pi_k(x_j) \in B(\ell^2(\cP_k))$ for $j=1,2,\ldots m$. 
Fix a subset $\cS_k$ of $\cP_k$ such that $|\cS_k| = \frac{|\cP_k|-1}{2}$.
Using the construction of Runde \cite{Runde}, we define $\pi_k(x_{m+1})\in B(\ell^2(\cP_k))$ by assigning
\[
\pi_k(x_{m+1})(e_\times) =
\begin{cases}
    e_\times , & e_\times\in \cS_k \\
    -e_\times, & e_\times\notin \cS_k
\end{cases}
\]
where $\{e_\times: \times \in \cP_k\}$ is the canonical basis of $\ell^2(\cP_k)$.

The next lemma, which is crucial for us, is from \cite{Runde}. It reveals an unusual property of $SL_3(\bZ)$.

Consider the following attempt for $SL_3(\bZ)$.
\begin{attempt}\label{L: task}
Given $\ep>0$, find $r\in \bN$ such that for each prime number $p_k$, one can get $\xi_1,\ldots,\xi_r, \eta_1,\ldots,\eta_r \in \ell^2(\cP_k)$ that satisfy
\[
\|\sum_{i=1}^r \xi_i\otimes \eta_i -(\pi_k(x_j)\otimes \pi_k(x_j))(\xi_i\otimes \eta_i)\|_{\ell^2(\cP_k)\hot^\gamma \ell^2(\cP_k)} \leq \ep \|\sum_{i=1}^r \xi_i\otimes \eta_i\|_{\ell^2(\cP_k)\hot^\gamma \ell^2(\cP_k)}
\]
for all $j=1,2,\ldots, m+1$.
\end{attempt}

\begin{lemma}[Lemma~3.4.5 of \cite{Runde}]\label{L: Runde}
    There is $\ep>0$ such that Attempt~\ref{L: task} cannot be attained.
\end{lemma}

\begin{prop}\label{P: not aa}
Let $\cM$ be a type I von Neumann algebra of the form \eqref{E: type I}.
Suppose that $\sup_{\al\in \Gamma}\dim(\cH_\al) =\infty$. Then $\cM$ is not approximately amenable.
\end{prop}

\proof

We may assume $\dim(\cH_\al)< \infty$ for all $\al$. Then we may write
\[
\cM = \oplus_{k=1}^\infty A_{n_k}\otimes B(\cH_{n_k})
\]
where each $A_k$ is a non-trivial commutative von Neumann algebra and $(n_k)_{k=1}^\infty$ is a strictly increasing sequence of positive integers. If $\cM$ is approximately amenable, then, as a homomorphic image, $\oplus_{k=1}^\infty B(\cH_{n_k})$ is approximately amenable. So we may assume 
\[
\cM = \oplus_{k=1}^\infty B(\cH_{n_k})
\]
and prove that it is not approximately amenable. 

Let $(p_k)_{k=1}^\infty$ be the increasing sequence of all prime numbers. Going through an isomorphic image, we can assume $n_k \geq |\cP_k|$ for all $k$.
Let $m_k = n_k - |\cP_k|$. Denote $\bN_n =\{1,2, \ldots n\}$ for $n\in \bN$ and $\bN_0 =\emptyset$. Then we can write 
\[\cH_{n_k} = \ell^2(\cP_k) \oplus_2 \ell^2(\bN_{m_k})
\]
and
\[
\cM = \oplus_{k=1}^\infty B\left(\ell^2(\cP_k) \oplus_2 \ell^2(\bN_{m_k})\right).
\]
Let
\[
\cH =\ell^2\text{-}\oplus_{k=1}^\infty \cH_{n_k} = 
\ell^2\text{-}\oplus_{k=1}^\infty\left(\ell^2(\cP_k) \oplus_2 \ell^2(\bN_{m_k})\right),
\]
and denote the canonical basis of $\cH$ by $E$.
To disprove the approximate amenability of $\cM$, we use the same argument as \cite[Theorem~3.4.12]{Runde}.

Consider the unitary representation $\omega_k=\pi_k\oplus \text{id}_k$ of $G$ on $\ell^2(\cP_k) \oplus_2 \ell^2(\bN_{m_k})$, where id$_k$ is the identity representation on $\ell^2(\bN_{m_k})$. 
Let $\pi_k(x_j)$, $j=1,2,\ldots m+1$, be the elements of $B(\ell^2(\cP_k))$ as defined before Attempt~~\ref{L: task}. 
We have 
\[
\Omega(x_j) =\oplus_{k=1}^{\infty} \omega_k(x_j) \in \cM
\]
for each $j=1,2,\ldots m+1$.

Suppose, to the contrary, that $\cM$ is approximately amenable. Then, for every given $\ep >0$, there is $T = \sum_{i=1}^r a_i\otimes b_i\in \cM \otimes\cM$ such that $\sum_{i=1}^r a_ib_i = \rm{id}$ and 
\[
\left\|\Omega(x_j) \cdot T -T \cdot \Omega(x_j)\right\|_\gamma <\frac{\ep}{m+1}
\]
for all $j  =1, 2, \ldots, m+1$.

Let $P_k$ be the orthogonal projection of $\cH$ on $\ell^2(\cP_k)$, and let $T_k = \sum_{i=1}^rP_k \cdot a_i\otimes P_k\cdot b_i^*$.
By Lemma~\ref{L: Oz}, for each $e\in E$,
  \[
   T_k(e) := \sum_{i=1}^rP_k \cdot a_i(e)\otimes P_k\cdot b_i^*(e)\in \ell^2(\cP_k)\hot^\gamma \ell^2(\cP_k)
    \]
and
    \[
    \sum_{e\in E} \|T_k(e)\|_{\ell^2(\cP_k)\hot^\gamma \ell^2(\cP_k)} \leq |\cP_k| \|T\|_\gamma.
    \]
We note that for each $x_j$,
\begin{align*}
    (\pi_k(x_j)\otimes \pi_k(x_j))T_k(e) &= \sum_{i=1}^r\pi_k(x_j)P_k  a_i(e)\otimes \pi_k(x_j)P_k b_i^*(e)\\
    & = \sum_{i=1}^rP_k [\Omega(x_j)
    \cdot  a_i](e)\otimes P_k[\Omega(x_j)
    \cdot b_i^*](e)\\
& = \sum_{i=1}^rP_k [\Omega(x_j)
\cdot  a_i]
(e)\otimes P_k 
[b_i\cdot \Omega(x_j)
]^*(e).
\end{align*}

Apply Lemma~\ref{L: Oz} again. We have
\begin{align*}
    \sum_{e\in E} \| &[T_k -(\pi_k(x_j)\otimes \pi_k(x_j)T_k](e)\|_{\ell^2(\cP_k)\hot^\gamma \ell^2(\cP_k)} \\
    &\leq |\cP_k| \|T - \Omega(x_j)   
    \cdot T \cdot \Omega(x_j)    
    ^{-1}\|_\gamma \\
    &\leq |\cP_k| \|\Omega(x_j)   
    \cdot T - T\cdot \Omega(x_j)    
    \|_\gamma < |\cP_k|\frac{\ep}{m+1}
\end{align*}
for all $j=1,2,\ldots, m+1$. Thus,
\[
 \sum_{e\in E} \sum_{j=1}^r \|[T_k -(\pi_k(x_j)\otimes \pi_k(x_j)T_k](e)\|_{\ell^2(\cP_k)\hot^\gamma \ell^2(\cP_k)} <|\cP_k|\ep .
 \]

On the other hand,
\begin{align*}
\sum_{e\in E} \|T_k(e)\|_{\ell^2(\cP_k)\hot^\gamma \ell^2(\cP_k)}
&\geq \sum_{e\in E} \sum_{i=1}^r \langle P_k \cdot a_i(e), P_k\cdot b_i^*(e)\rangle \\
& = \sum_{i=1}^r {\rm tr} (b_i \cdot P_k \cdot a_i) = {\rm tr} ((\sum_{i=1}^r a_ib_i )P_k) = {\rm tr} (P_k) = |\cP_k|.
\end{align*}
Therefore,
\begin{align}\label{E: ineq for sum}
\sum_{e\in E}\sum_{j=1}^{m+1} &\| [T_k -(\pi_k(x_j)\otimes \pi_k(x_j)T_k](e)\|_{\ell^2(\cP_k)\hot^\gamma \ell^2(\cP_k)} \notag \\
&\leq \ep 
\sum_{e\in E} \|T_k(e)\|_{\ell^2(\cP_k)\hot^\gamma \ell^2(\cP_k)}
\end{align}
This is true for all $k\in \bN$. So, for each $k$, there is $e_k\in E$ for which
\begin{equation}\label{E: ineq}
    \sum_{j=1}^{m+1} \| T_k(e_k) -(\pi_k(x_j)\otimes \pi_k(x_j)T_k(e_k)\|_{\ell^2(\cP_k)\hot^\gamma \ell^2(\cP_k)}
\leq \ep 
 \|T_k(e_k)\|_{\ell^2(\cP_k)\hot^\gamma \ell^2(\cP_k)}.
\end{equation}

 If the right side is $0$, then so is the left side. We may remove the corresponding terms from the original relation \eqref{E: ineq for sum} and reselect $e_k$ so that \eqref{E: ineq} holds with $T_k(e_k)\neq 0$.
This ensures that, for each $p_k$, we may choose $\sum_{i=1}^r \xi_i\otimes \eta_i = T_k(e_k) \neq 0$ so that Attempt~\ref{L: task} can be achieved. This contradicts Lemma~\ref{L: Runde}, since $\ep>0$ was arbitrary.  The proof is complete.
\qed

Combining Propositions~\ref{P: dim bdd} and \ref{P: not aa}, we derive the following result.

\begin{thm}
A type I von Neumann algebra is approximately amenable if and only if it is amenable.
\end{thm}

\begin{remark}
The question of whether $\ell^\infty (\bM_n)$ is approximately amenable was raised in \cite[Page 250]{G-L04}. We note that this algebra is indeed the type I von Neumann $\ell^\infty \text{-} \oplus_{n=1}^\infty B(H_n)$, where $H_n$ denotes the $n$-dimensional Hilbert space. Theorem~\ref{P: not aa} answers this question negatively.
\end{remark}


\begin{thebibliography}{9}
\bibitem{G-L04}
F. Ghahramani and R. J. Loy, Genaralized notions of amenability,
JFA 208 (2004),
229-260.

\bibitem{GLZ08}
  F. Ghahramani, R. J. Loy and Y. Zhang,
Generalized notions of amenability, II,  J. Funct. Anal. 254 (2008), 1776-1810.

\bibitem{GZ07}
F. Ghahramani and Y. Zhang, 
Pseudo-amenable and pseudo-contractible Banach algebras. Math. Proc. Cambridge Phil. Soc. 142 (2007), 111-123.

\bibitem{Oz}
N. Ozawa,
A note on non-amenability of $B(\ell_p)$ for $p = 1,2$, Internat. J. Math. 15 (2004), 557-565.
\bibitem{Runde}
V. Runde,
Amenable Banach algebras, Springer, 2020.

\bibitem{Tak I}
M. Takesaki,
Theory of Operator Algebras I,
Springer, 1979.
\bibitem{Z22}
Y. Zhang, 
Approximate amenability and pseudo-amenability in Banach algebras,
AMSA 8 (2023), 309-320.
\end{thebibliography}
\end{document}